\let\beginformat=\relax
\let\endformat=\relax

\batchmode

\beginformat
\documentclass[10pt,reqno,oneside]{amsart}
\usepackage{amssymb}
\usepackage{amsmath}
\usepackage{amsthm}
\usepackage{graphicx}
\endformat

\theoremstyle{plain}
\newtheorem{Thm}{Theorem}
\newtheorem{Cor}[Thm]{Corollary}

\newtheorem{\theMain}{}
\newtheorem{Lem}[Thm]{Lemma}
\newtheorem{Prob}{Problem}
\newtheorem{Prop}[Thm]{Proposition}

\theoremstyle{definition}
\newtheorem{Def}{Definition}
\newtheorem*{note}{Note}
\def\determinant#1{\left|#1\right|}
\def\biggdet#1{\biggl|#1\biggr|}
\def\card#1{\left|#1\right|}

\begin{document}

%Topmatter
\title[Enumeration of Tilings of Diamonds and Hexagons]{Enumeration of 
Tilings of Diamonds and Hexagons \\ with Defects}
\author{Harald Helfgott and Ira M. Gessel}
\address{Mathematics Department \\
  Brandeis University
  Waltham, MA 02254-9110}
\email{hhelf@cs.brandeis.edu, gessel@math.brandeis.edu}
\thanks{Harald Helfgott's research was funded in part by the Undergraduate
Research Program of Brandeis University. Ira Gessel's research was supported
by NSF grant DMS9622456.}
\date{July 30, 1998}

%End topmatter

\begin{abstract} 
We show  how to count tilings of Aztec diamonds and hexagons
with defects using determinants. In several cases these determinants can be
evaluated in closed form.  In particular, we obtain solutions to  open
problems 1, 2, and 10 in James Propp's list of problems on enumeration of
matchings~\cite{jP1}. 
\end{abstract}
\maketitle
\section{Introduction}

%\begin{center}
%\begin{minipage}[b]{0.4\textwidth}
%\centering \includegraphics{aztec3.eps}
%\par\vspace{0pt}
%\end{minipage}
%\begin{minipage}[b]{0.4\textwidth}
%\centering \includegraphics[angle=90]{hexa34.eps}
%\par\vspace{0pt}
%\end{minipage}
%\caption{Aztec Diamond of Order 3} \label{fig:side:a}
%\end{center}

While studying dimer models, P. W. Kasteleyn  \cite{K} noticed that tilings of
very simple figures by very simple tiles can be not only plausible
physical models, but also  starting points for some very interesting
enumeration problems. Kasteleyn himself solved the problem of counting
tilings of a rectangle by dominoes. He also found a general method
(now known as Kasteleyn matrices) for computing the number of tilings of
any bipartite planar graph in polynomial time. Kasteleyn's method has
proven very useful for computational-experimental work, but it does
not, of itself, provide proofs of closed formulas for specific
enumeration problems.  We shall see a few examples of problems for
which Kasteleyn matrices alone are inadequate.

\begin{figure}
        \begin{minipage}[b]{0.5\linewidth}
                \centering \includegraphics{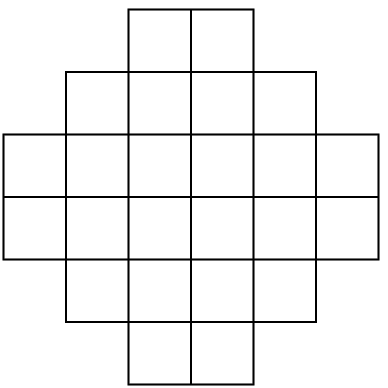}
                \caption{\hbox{Aztec diamond of order 3}}
\label{fig:side:a}   
        \end{minipage}%
\hglue-30pt
        \begin{minipage}[b]{0.5\linewidth}
                \centering \includegraphics[height=1.5in]{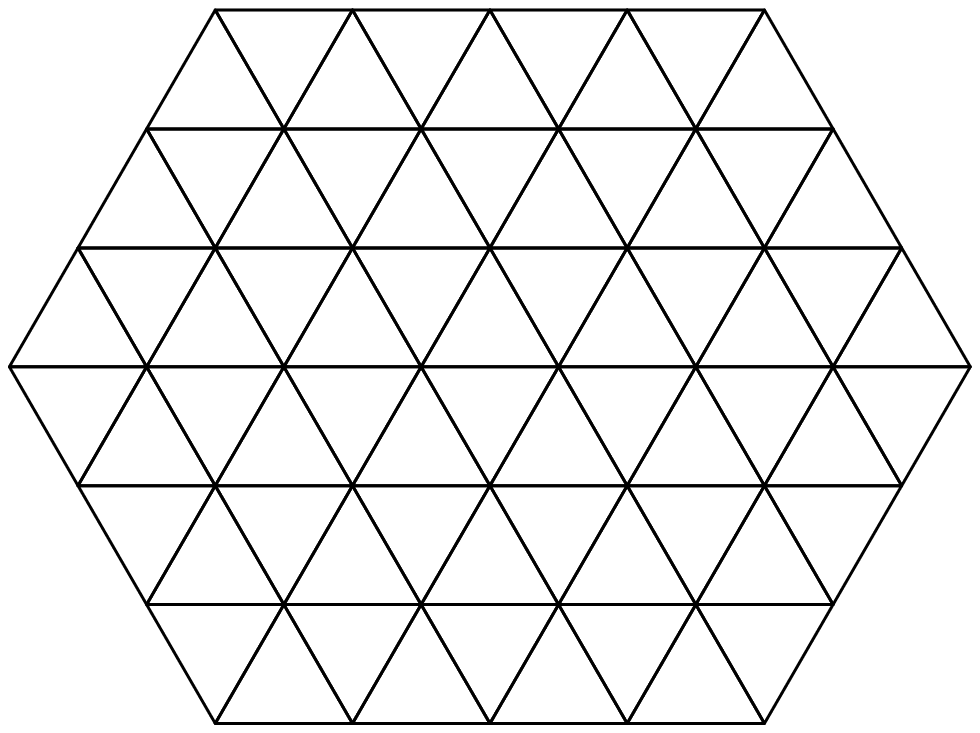}
                \caption{$(3,4,3)$ hexagon} \label{fig:side:b}
        \end{minipage}%
\end{figure}

By an {\it $(a,b,c,d,e,f)$ hexagon} we mean a hexagon with sides of
lengths $a,b,c,d,e,f$, and angles of 120 degrees, subdivided into
equilateral triangles of unit side by lines parallel to the sides.  We
draw such a hexagon with the sides of lengths $a,b,c,d,e,f$ in
clockwise order, so that the side of length $b$ is at the top and the
side of length $e$ is at the bottom. We shall use the term {\it
$(a,b,c)$ hexagon \/} for an $(a,b,c,a,b,c)$ hexagon. Thus
Figure~\ref{fig:side:b} shows a $(3,4,3)$ hexagon.

An {\it Aztec diamond of order $n$} is the union of all unit squares
with integral vertices contained within the region $|x|+|y| \leq
n+1$. Figure ~\ref{fig:side:a} shows an Aztec diamond of order 3.

We are interested in tilings of hexagons with {\it lozenges\/}, which are
rhombi with unit sides and angles of 120 and 60 degrees, and tilings of Aztec
diamonds with {\it dominoes\/}, which are  1 by 2 rectangles. In particular,
we shall examine three problems from James Propp's list of open problems on
tilings \cite{jP1}.

\begin{Prob}[Propp's Problem 1]\label{Propp1}
Show that in the $(2n-1,2n,2n-1)$ hexagon, the central
vertical lozenge (consisting of the two innermost triangles) 
is covered by a lozenge in exactly one-third of the tilings.
\end{Prob}

\begin{Prob}[Propp's Problem 2]\label{Propp2}
Enumerate the lozenge-tilings of the region obtained from the
$(n, n+1, n, n+1, n, n+1)$ hexagon by removing the central triangle.
\end{Prob}

\begin{Prob}[Propp's Problem 10]\label{Propp10}
Find the number of domino tilings of a $(2k-1)$ by $2k$
undented Aztec rectangle with a square adjoining the central square removed,
where the $a$ by $b$ undented Aztec rectangle is defined as the union of the
squares bounded by $x+y\leq b+1$, $x+y\geq b-2a-1$, $y-x\leq b+1$,
$y-x\geq -(b+1)$.
\end{Prob}

We have solved these three problems, not by using Kasteleyn matrices,
but by choosing a new approach, which, while much less general than
Kasteleyn matrices, is better suited for problems like these three.
We can summarize our approach as follows:

\begin{enumerate}
\item Find the number of tilings of half of a hexagon or half of
a diamond, with dents at given places. This is not new: 
see \cite{CLP} and \cite{EKLP}. 
\item Express the number of tilings of the figure as a whole as a sum
of squares of the expressions obtained in the first step. The sum's range
depends on the ``defects'' (missing triangles or squares, fixed
lozenges or dominos) given in the problem.
\item Express the sum of squares as a Hankel determinant.
\item Evaluate the Hankel determinant using continued fractions or
Jacobi's theorem. 
\end{enumerate}

C. Krattenthaler has been working on these problems at the same time
as us, together with M. Ciucu \cite{CK} and S. Okada \cite{KO}. The
solution to Problem 1 in \cite{CK} is literally orthogonal to ours:
Ciucu and Krattenthaler slice the hexagon vertically rather than
horizontally. More generally, Fulmek and Krattenthaler \cite {FK}
have counted tilings of an $(n,m,n,n,m,n)$ hexagon  that contain an
arbitrary fixed rhombus on the symmetry axis that cuts through the
sides of length $m$.  Krattenthaler and Okada's solution \cite{KO} to
Problem 2 and Krattenthaler's solution \cite{Kr} to Problem 10  are much like
ours in steps 1 and 2. Thereafter, they are based on identities for Schur
functions, not Hankel determinants.  The work of Krattenthaler and his coauthors  and our
work thus complement each other.

\section{From Tilings to Determinants}

First we note that a necessary and sufficient condition for an
$(a,b,c,d,e,f)$ hexagon to exist is that the parameters be
nonnegative integers satisfying
$a-d=c-f=e-b$.  The number of upward pointing triangles  minus the
number of downward pointing triangles in an $(a,b,c,d,e,f)$
hexagon is $a-d$. Then since every lozenge covers one upward
pointing triangle and one downward pointing triangle, an
$(a,b,c,d,e,f)$ hexagon can be tiled by lozenges only if $a=d$,
and this implies that that the hexagon is an $(a,b,c)$ hexagon.
Moreover, if we remove $a-d$ upward pointing triangles from an
$(a,b,c,d,e,f)$ hexagon with $a\ge d$, then the remaining figure
will have as many upward pointing as downward pointing triangles.

\begin{figure}
        \begin{minipage}[b]{0.5\linewidth}
                \centering \includegraphics{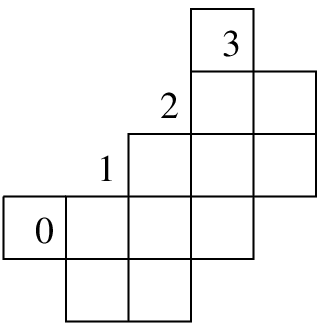}
                \caption{Dented 3 by 2 Aztec rectangle} \label{fig:side:c}      
        \end{minipage}%
  \hglue -30pt
        \begin{minipage}[b]{0.5\linewidth}
         \centering \includegraphics[height=0.8in]{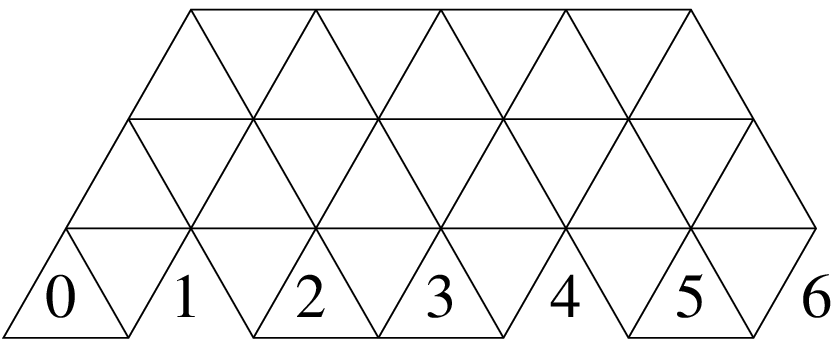}
                \caption{Dented $(3,4,3)$ semi-hexagon}
\label{fig:side:d}
        \end{minipage}%
\end{figure}

\begin{Def}
        A {\it $(k,q,k)$ upper semi-hexagon\/} 
is the upper half  of a 
$(k,q,k)$ hex\-a\-gon
having sides $k,q,k,q+k$, i.e., a symmetric trapezium. A {\it $(k,q,k)$ lower
semi-hexagon\/} is defined similarly. A {\it $(k,q,k)$ dented upper semi-hexagon\/}
is a $(k,q,k)$ semi-hexagon with $k$ upward pointing triangles removed from the 
side of length $q+k$. (Figure~\ref{fig:side:d} shows a $(3,4,3)$ dented upper
semi-hexagon with dents at positions 1, 4, and 6). It will be convenient to use
the term {\it semi-hexagon\/} for an upper semi-hexagon.
\end{Def}

Note that a $(k,q,k)$ semi-hexagon is the same as a $(k,q,k,0,q+k,0)$ hexagon, so removing 
$k$ upward pointing triangles leaves a region with as many upward as downward triangles.

\begin{Def}
        An $a$ by $b$ dented Aztec rectangle is the union of the
squares bounded by $x+y\leq b+1$, $x+y\geq b-2a-1$, $y-x\leq b$,
$y-x\geq -(b+1)$, with the squares in positions
$r_0<r_1<\cdots< r_{b-1}$ removed from the side given by $y-x\leq b$
(see Figure~\ref{fig:side:c}).
\end{Def}

Before proceeding with our results on tilings, we first note some facts about the
power sums $1^j+2^j+\cdots+m^j$ that we will need later on. We omit the
straightforward proofs.

For any integer $m$ and any nonnegative integer $j$ we define $S_m^j$ by 
\[ 
S_m^j=\begin{cases}
   1^j+\cdots+m^j, &\text{if $m>0$;}\\
   0, &\text{if $m=0$;}\\
   (-1)^{j+1}\bigl(0^j+1^j+\cdots+(-m-1)^j\bigr), &\text{if $m<0$,}
\end{cases}
\]
where we interpret $0^0$ as 1.

\renewcommand{\labelenumi}{(\theenumi)}
\begin{Lem}\label{Lem:0}
The numbers $S_m^j$ have the following properties:
\begin{enumerate}
\item For any integers $p$ and $q$, with $p\le q$, 
  \[p^j+(p+1)^j+\cdots+q^j=S_q^j-S_{p-1}^j.\]
\item $S_0^j=0$ for all $j$ and $S_{-1}^j=0$ for $j>0$.
\item  For $m>0$, $S_{-m}^j = (-1)^{j+1}S_{m-1}^j$.
\item For $m\ge0$, $S_m^j$ is given by the exponential generating function
\[\sum_{j=0}^\infty S_m^j \frac {x^j}{j!} =\frac{e^x(e^{mx}-1)}{e^x-1}\]
\item $ S_m^j$ is a polynomial in $m$ of degree $j+1$, with leading coefficient $1/(j+1)$.
\end{enumerate}
\end{Lem}

Next we prove two known results. First, we have
a closed expression for the number of tilings of semi-hexagons 
with given dents, first stated in this form in~\cite{CLP}. This is equivalent to a
well-known result on the enumeration of Gelfand patterns, as noted in \cite{CLP}, or
on column-strict plane partitions. (See Knuth
\cite[exercise 23, p.~71; solution, p.~593]{Kn} for a proof similar to ours.)

\begin{Lem}\label{Lem:1}
        The number of tilings of a $(k,q,k)$ semi-hexagon
with dents at positions $0\le
r_0<\cdots<r_{k-1}< q+k$ is
\[T_{k,q,r} = \frac{1}{V_{k-1}}  \prod_{0\leq i<j<k} (r_j-r_i),\]
where $V_n=1!\, 2!\cdots n!=\prod_{0<i<j\le n}(j-i)$.
\end{Lem}

\begin{proof}
        We proceed by induction on $k$. For the case $k=1$,
there is only one tiling, no matter where
the solitary dent is.  Hence the lemma holds for $k=1$.

Let us now assume the lemma holds for $k$. Suppose we have
a tiling of a $(k+1,q,k+1)$  with dents at
$0\leq r_0< \dots < r_k< q+k+1$. If we remove the bottom layer
of lozenges from the dented side, we obtain a tiling of a
$(k,q,k)$ semi-hexagon with dents at $0\leq t_0< \dots < t_{k-1}\leq q+k$,
$r_i\leq t_i < r_{i+1}$. For every such tiling of a $(k,q,k)$
semi-hexagon
with dents at those places, there is exactly one tiling of the
dented $(k,q,k)$ semi-hexagon. Hence
\begin{align*}
T_{k+1,q,r} &= \sum_{r_i\leq t_i< r_{i+1}} T_{k,q,t} \\
&= \sum_{r_i\leq t_i < r_{i+1}} \frac{1}{V_{k-1}}  
    \prod_{0\leq i<j<k} (t_j-t_i). \\
&= \frac{1}{V_{k-1}}
   \sum_{r_i\leq t_i < r_{i+1}} \left| t_i^j \right|_0^{k-1} \\
&= \frac{1}{V_{k-1}}
   \left| S_{r_{i+1}-1}^j-S_{r_i-1}^j \right|_0^{k-1} \\
&= \frac{1}{V_{k-1}}
   \left| S_{r_{i+1}-1}^j-S_{r_0-1}^j \right|_0^{k-1},
\end{align*}
where $S_m^j=1^j + 2^j+\dots + m^j$. In the second line of our
calculations we can see that, since $\prod_{0\leq i<j<k} (t_j-t_i) $
depends only on the differences between the $t_i$'s, $T_{k+1,q,r}$
depends only on the differences between the $r_i$'s, not on their
actual values. (It is also easy to see this combinatorially.) Hence it is sufficient to prove
the formula in the case
$r_0=0$. By Lemma \ref{Lem:0}, $S_{m-1}^j-S_{-1}^j$ is a polynomial in $m$ of degree $j+1$
with leading coefficient $1/(j+1)$ that vanishes at~$0$. Thus we can reduce the determinant  
$\left|
S_{r_{i+1}-1}^j-S_{-1}^j \right|_0^{k-1} $ to $\left| r_{i+1}^{j+1}/(j+1) \right|_0^{k-1}$ by
elementary column operations. Hence
$$
\begin{aligned}
T_{k+1,q,r} &= \frac{1}{V_{k-1}} 
\left| \frac{r_{i+1}^{j+1}}{j+1} \right|_0^{k-1} \\
            &= \frac{1}{V_k}
\left| r_{i+1}^{j+1} \right|_0^{k-1} \\
            &= \frac{r_1  r_2 \dotsb r_k}{V_k} 
\left| r_{i+1}^j \right|_0^{k-1} \\
            &= \frac{r_1  r_2 \dotsb r_k}{V_k} 
\prod_{1\leq i<j<k+1} (r_j-r_i) \\
&=\frac1{V_k}\prod_{0\leq i<j<k+1} (r_j-r_i),
\end{aligned}
$$
since we assumed that $r_0=0$. Then by our observation the formula holds for all values of
$r_0$.
\end{proof}

\section{Tilings of Dented Aztec Rectangles}

\begin{Def}
        An $a$ by $b$ dented Aztec rectangle is the union of the
squares bounded by $x+y\leq b+1$, $x+y\geq b-2a-1$, $y-x\leq b$,
$y-x\geq -(b+1)$, with the squares in positions
$r_0<r_1<\cdots< r_{b-1}$ removed from the side given by $y-x\leq b$
(see Figure~\ref{fig:side:c}). An $a$ by $b$ undented Aztec rectangle 
is an $a$ by $b+1$ dented Aztec rectangle with all squares
on the side given by $y-x\leq b$ removed.
\end{Def}

Our next result counts dented Aztec rectangles. Another proof can be found
in \cite{EKLP}. Just as tilings of dented hexagons correspond to Gelfand
patterns, in \cite{EKLP} it is shown that tilings of  dented Aztec
rectangles correspond to monotone triangles, and in this context, a proof
of the formula can be found in \cite{MRR}.

\begin{Lem}\label{Lem:3}
        The number of tilings of an $a$ by $b$ dented Aztec rectangle  
with dents at $0\leq r_0\le \dots \le r_{b-1}\leq a$ is
$$A_{a,b,r} = \frac{2^{\frac{b(b-1)}{2} }}{V_{b-1}}
 \prod_{0\leq i<j<b} (r_j-r_i),$$
where $V_n=1!\, 2!\cdots n!$.
\end{Lem}
\begin{proof}
We proceed by induction on $b$. 
First we note that if $r_i=r_{i+1}$ for some $i$, then
the lemma asserts that
$A_{a,b,r}=0$, which is correct. Although of no interest in itself, this case will
be necessary for the induction.

If
$b=1$, there is only one tiling, no matter where the one dent is.
(In general, the number of dents
has to be equal to $b$ for the dented Aztec rectangle to be tileable.)
Hence the lemma holds for $b=1$.

Let us now assume the lemma holds for $b$.  
Suppose we have
a tiling of an $a$ by $b+1$ Aztec rectangle with dents at
$0\leq r_0< \dots < r_b\leq a$. 
If we remove all dominoes with one or two squares on the dented long
diagonal and the adjacent short diagonal,
we obtain a tiling of an
$a$ by $b$ Aztec rectangle with dents at $0\leq t_0< \dots < t_{b-1}\leq a$,
where $r_k\leq t_k \leq r_{k+1}$. For every such tiling of an $a$ by $b$
Aztec rectangle with dents at those places, there are $2^m$ tilings
of the $a$~by $b+1$ dented Aztec rectangle, where $m$ is the
cardinality of
$\{ k: r_k < t_k < r_{k+1} \}$. 

Next we show that this implies  
\begin{equation}\label{Eq:01}
A_{a,b+1,r} = \sum_{l\in \{0,1\}^b} 
                \sum_{r_k\leq t_k-l_k< r_{k+1}} A_{a,b,t}
\end{equation}
This follows from the fact that if $r_k<t_k< r_{k+1}$ then $r_k\leq t_k-l_k<
 r_{k+1}$  if $l_k$ is either 0 or 1, but if $r_k=t_k$ then  this
inequality holds only for $l_k=0$  and if $r_k= t_{k+1}$, it holds only for $l_k=1$.
Thus the number of different possible values of $l$ corresponding to a given sequence
$r$, is $2^m$, where $m$ is the
cardinality of $\{ k: r_k < t_k < r_{k+1} \}$. 
Moreover, if for some $l\in \{0,1\}^b$, $t$ satisfies $r_k\leq t_k-l_k< r_{k+1}$
for all $i$, then we must have $t_0\le t_1 \le\cdots\le t_{b-1}$, so
all terms $A_{a,b,t}$ that occur  in \ref{Eq:01} either have $t_0< \cdots <
t_{b-1}$ or are zero; in either case they are covered by the
induction hypothesis.

Hence
\begin{align}
A_{a,b+1,r} &= \sum_{l\in \{0,1\}^b} 
                 \sum_{r_k\leq t_k-l_k< r_{k+1}} A_{a,b,t}\notag \\  
            &= \frac{2^{\frac{b(b-1)}{2}}}{V_{b-1}} 
                \sum_{l\in \{0,1\}^{b}} 
                \sum_{r_k+l_k\leq t_k< r_{k+1}+l_k} 
                \prod_{0\leq i<j<b} (t_j-t_i) \label{Eq:tdiff}\\
            &= \frac{2^{\frac{b(b-1)}{2}}}{V_{b-1}}
                \sum_{l\in \{0,1\}^{b}} 
                \sum_{r_k+l_k\leq t_k< r_{k+1}+l_k} 
                \left| t_i^j \right|_0^{b-1} \notag\\
            &= \frac{2^{\frac{b(b-1)}{2}}}{V_{b-1}}
                \sum_{l\in \{0,1\}^{b}} 
                \left| S_{r_{i+1}+l_{i}-1}^j - S_{r_i+l_i-1}^j \right|_0^{b-1}, \notag
\end{align}
where $S_m^j=1^j + 2^j+\dots + m^j$.
Now if $u(i,j,k)$ is any function defined for $0\le i,j < b$, $0\le k\le 1$, then  since a
determinant is a linear function of its rows, we have
\[\sum_{l\in \{0,1\}^{b}} \left| u(i,j,l_i)\right|_0^{b-1}
  = \left| u(i,j,0)+u(i,j,1)\right|_0^{b-1}.\]
Thus
\begin{align*}
 A_{a,b+1,r} &= \frac{2^{\frac{b(b-1)}{2}}}{V_{b-1}}
                \left| S_{r_{i+1}-1}^j + S_{r_{i+1}}^j 
                    - (S_{r_i-1}^j + S_{r_i}^j) \right|_0^{b-1} \\
            &= \frac{2^{\frac{b(b-1)}{2}}}{V_{b-1}}
            \left| (S_{r_{i+1}-1}^j + S_{r_{i+1}}^j) - 
                   (S_{r_0-1}^j + S_{r_0}^j) \right|_0^{b-1},
\end{align*}

 By (\ref{Eq:tdiff}), we can see that, since $\prod_{0\leq i<j<b} (t_j-t_i) $
depends only on the differences between the $t_k$'s, $A_{a,b+1,r}$
depends only on the differences between the $r_k$'s, not on their
actual values. Hence we may assume that $r_0=0$. 
Since $S_{m-1}^j + S_m^j-(S_{-1}^j+S_0^j)=S_{m-1}^j + S_m^j-S_{-1}^j$ is a polynomial
in $m$ of degree $j+1$ with leading coefficient $2/(j+1)$ that vanishes at~$0$, we can
reduce the determinant $\left| (S_{r_{i+1}-1}^j + S_{r_{i+1}}^j)-S_{-1}^j\right|_0^{b-1} $
to $\left| 2 r_{i+1}^{j+1}/(j+1) \right|_0^{b-1}$ by elementary column operations.
Hence
$$
\begin{aligned}
A_{a,b+1,r} &=  \frac{2^{\frac{b(b-1)}{2}}}{V_{b-1}}
                \left| \frac{2 r_{i+1}^{j+1}}{j+1} \right|_0^{b-1} \\
            &=  \frac{2^{\frac{(b+1)b}{2}}}{V_b}
                \left| r_{i+1}^{j+1} \right|_0^{b-1} \\
            &=  \frac{2^{\frac{(b+1)b}{2}}r_1  r_2\dotsb r_b}{V_b} 
                \left| r_{i+1}^j \right|_0^{b-1} \\
            &=  \frac{2^{\frac{(b+1)b}{2}} r_1  r_2\dotsb r_b}{V_b}
                \prod_{1\leq i < j < b+1} (r_j-r_i)\\
            &=  \frac{2^{\frac{(b+1)b}{2}}}{V_b}  
                \prod_{0\leq i < j < b+1} (r_j-r_i) .
\end{aligned}
$$
\end{proof}

\section{From hexagons to Determinants}

We now compute the number of tilings of a $(k,q,k)$ hexagon with 
restrictions on where vertical lozenges may cross the horizontal symmetry 
axis. 

\begin{Prop}
\label{Prop:L}
 Let $L$ be a subset of $\{0,1,\dots, k+q-1\}$. Then the number
of tilings of a $(k,q,k)$ hexagon in which
the set of indices of the vertical lozenges crossing the
$q+k$-long symmetry axis is a subset of $L$ is
\[\frac{1}{V_{k-1}^2}  \biggdet{ \sum_{l\in L} l^{i+j} }_0^{k-1}. \]
\end{Prop}
\begin{proof}
We first recall that by the Binet-Cauchy theorem \cite[p.~9]{G},
if  $M$ is any $k$ by $n$ matrix and $M^t$ is its transpose, then the determinant
of $MM^t$ is equal to the sum of the squares of the
$k$ by $k$ minors of $M$.

The number of tilings of a $(k,q,k)$ hexagon in which the
indices of the vertical lozenges crossing the
$q+k$-long symmetry axis are $r_0<r_1<\cdots<r_{k-1}$ is clearly 
\[T_{k,q,r}^2=\frac 1{V_{k-1}^2}\bigl(\determinant{r_j^i}_0^{k-1}\bigr)^2.\]
Thus the number of tilings to be counted is the sum of $T_{k,q,r}^2$ over all
$r_0<r_1<\cdots<r_{k-1}$ where each $r_j$ is in $L$. Now suppose that the
elements of $L$ are $l_0<l_1<\cdots<l_{n-1} $ and let $M$ be the $k$ by
$n$ matrix
$(l_j^i)_{0\le i<k,\, 0\le j < n}$. Then by the Binet-Cauchy theorem,
\[\sum_r T_{k,q,r}^2=\frac 1 {V_{k-1}^2}\determinant{M M^t}
  =\frac 1 {V_{k-1}^2}\biggdet{\sum_{l\in L}l^{i+j}}_{0}^{k-1}.\]
\end{proof}

Note that since the numbers $T_{k,q,r}$ depend only on the differences of the
$r_i$,  the determinant in Proposition \ref{Prop:L} depends only on the
differences of the elements of
$L$; thus we may shift all the elements of $L$ by the same amount without changing the
determinant. This observation will be useful later on:

\begin{Lem} \label{Lem:Shift}
For any finite set $L$ of numbers and any number $u$,
\[\determinant{\sum_{l\in L} l^{i+j}}_0^{k-1}
 =\determinant{\sum_{l\in L} (l+u)^{i+j}}_0^{k-1}.\]\qed
\end{Lem}

\begin{Prop}
\label{Prop:L2}
The number of tilings of a $(k,2n+1-k,k,k+1,2n-k,k+1)$ hexagon with
a triangle removed below the center of the horizontal line
dividing the two ``hemispheres'' is
$$\frac{1}{V_{k-1}  V_k} 
\biggdet{(1+(-1)^{i+j}) \biggl(\sum_{l=1}^{n} l^{i+j+1}\biggr) }_0^{k-1}.$$
\end{Prop}
\begin{proof}
If we cut such a tiled hexagon into two parts by a horizontal line through the
middle 
vertices, and then remove the lozenges that are bisected by this line, we 
obtain
a tiling of a $(k,2n+1-k,k)$ upper semi-hexagon with dents at points
$r_0<r_1<\cdots<r_{k-1},$ and a tiling of a $(k+1,2n-k,k+1)$
lower semi-hexagon  with dents at points $r_0<r_1<\cdots<r_{k-1}$ and at the
center. Since the formula in Lemma~\ref{Lem:1} depends only on the
differences among the $r_i$, we can make zero lie on the center of
horizontal line dividing the two ``hemispheres'' of the hexagon.
Thus, we have $-n\leq r_0<r_1<\cdots<r_{k-1}\leq n,$ $r_i\ne 0.$

The number of tilings of the upper semi-hexagon is
$$\frac{1}{V_{k-1}} \prod_{0\leq i<j<k} (r_j-r_i) ,$$
and the number of tilings of the lower semi-hexagon is
$$\frac{1}{V_k} \prod_{0\leq i<j<k} (r_j-r_i) \prod_{0\leq i<k} |r_i|.$$
Hence the number of tilings of the hexagon for given 
$-n\leq r_0<\cdots<r_{k-1}\leq n$ is
$$\frac{1}{V_{k-1} V_k} (\left| r_j^i \right|_0^{k-1})^2 
                         \prod_{0\leq i<k} |r_i|
= \frac{1}{V_{k-1} V_k} (\left| |r_j|^{1/2} r_j^i \right|_0^{k-1})^2.$$ 
(Note that this vanishes whenever $r_j=0$ for some $j$.)

Now let $M$ be the $k$ by $2n+1$ matrix 
$(|j|^{\frac{1}{2}} j^i)_{0\leq i<k, -n\leq j \leq n}.$ 
Then by the Binet-Cauchy
theorem, the number of tilings of the hexagon is
$$
\begin{aligned}
\sum_{-n\leq r_0<\cdots<r_{n-1}\leq n}
\frac{1}{V_{k-1} V_k} \left(\left| |r_j|^{1/2} r_j^i \right|_0^{k-1}\right)^2
&= \frac{1}{V_{k-1} V_k} \left| M M^t \right| \\
&= \frac{1}{V_{k-1} V_k} 
\left| \sum_{-n\leq l \leq n} |l| l^{i+j} \right|_0^{k-1} \\
&= \frac{1}{V_{k-1} V_k} 
\left| (1+(-1)^{i+j}) \sum_{l=1}^n l^{i+j+1} \right|_0^{k-1} .
\end{aligned}
$$
\end{proof}

\section{From Aztec rectangles to determinants}

For our next result, we use the following lemma, which is analogous to the 
Binet-Cauchy theorem.

\begin{Lem}\label{Lem:PS}
Let $U=(u_{ij})$ be a $2k$ by $k$ matrix, with rows indexed from $0$ to $2k-1$ and columns
from $0$ to $k-1$. For each $k$-subset $A$ of $\{0,1,\ldots, 2k-1\}$, let $U_A$ be the $k$ by
$k$ minor of $U$ corresponding to the rows in $A$ and all columns, and let $\bar A$ be the
complement of $A$ in $\{0,1,\ldots, 2k-1\}$. Then
\[\sum_{\substack{A\subseteq \{0,\ldots,2k-1\}\\ |A|=k}} U_A U_{\bar A}
=2^k\determinant{u_{2i,j}}_0^{k-1} \determinant{u_{2i+1,j}}_0^{k-1}.
\]
\end{Lem}

\begin{proof} The lemma is a direct consequence of a result of Propp and Stanley
\cite[Theorem 2]{PS}. More precisely, the lemma follows from 
their result when we sum over all possibilities for $A^*$.
(As noted by Propp and Stanley, their result is a special case of a theorem of Sylvester
\cite{Sy}.)
\end{proof}

\begin{Prop}
        The number of tilings of an $a$ by $b$ undented Aztec rectangle, where
$a<b\leq 2a+1$, and  $b=2k+1$, with squares with indices $r_0<r_1<\dots <r_{b-a-1}$
missing from a diagonal of length $a+1$ 
going through the central square, is
\[
\frac{2^{k^2+a}}{V_k^2}
\biggl(
\prod_{0\leq j<i<2k+1-a} (r_i-r_j) 
\prod_{\substack{0\leq i<2a-2k\\[1pt] 0\leq j<2k+1-a}} \left| t_i-r_j \right|
\biggr)
\determinant{ t_{2i}^j }_0^{a-k-1} 
\determinant{ t_{2i+1}^j }_0^{a-k-1},
\]
where $t_0<t_1<\dots <t_{2a-2k-1}$ are the elements of $\{ 0,1,\dotsb a\} -
\{r_0,r_1,\dots ,r_{2k-a} \}$.
\end{Prop}
\begin{proof}
        Every tiling of the undented Aztec
rectangle with missing squares can be subdivided into a tiling
of two $a$ by $k+1$ dented Aztec rectangles with sets of dents disjoint
from each other and from $\{r_0,r_1,\dots ,r_{2k-a} \}$.
Let $T=\{ 0,1,\dotsb a\} -
\{r_0,r_1,\dots ,r_{2k-a} \}$, and let $t_0<t_1<\dots <t_{2a-2k-1}$ be the
elements of $T$. Then the number of tilings of the undented
Aztec rectangle with missing squares is
\begin{multline*}
\qquad\frac{2^{k(k+1)}}{V_k^2}
\prod_{0\leq j<i<2k+1-a} (r_i-r_j ) 
  \prod_{\substack{0\leq i<2a-2k\\[1pt] 0\leq j<2k+1-a}} 
  \left| t_i-r_j \right| \\
\times
 \sum_{\substack {P\subseteq T\\ \card {P}=a-k}}
\prod_{\substack{t_0,t_1\in P\\t_0<t_1}} (t_1-t_0)
\prod_{\substack{t_0,t_1\in T-P\\t_0<t_1}} (t_1-t_0),\qquad
\end{multline*}
which, written with determinants instead of products, is
$$
\frac{2^{k(k+1)}}{V_k^2} 
\biggl(
\prod_{0\leq j<i<2k+1-a} (r_i-r_j ) 
\prod_{\substack{0\leq i<2a-2k\\[1pt] 0\leq j<2k+1-a}} \left| t_i-r_j \right| \biggr)
\sum_{\substack{P\subseteq T\\ \card{P}=a-k}}
\determinant{p_i^j}_0^{a-k-1}
\determinant{q_i^j}_0^{a-k-1},
$$
where $p_0<p_1<\dots < p_{a-k-1}$ are the elements of $P$ and
$q_0<q_1<\dots < q_{a-k-1}$ are the elements of $T-P$.
Applying Lemma \ref{Lem:PS} yields the theorem.
\end{proof} 

We can prove the following proposition in exactly the same way.

\begin{Prop}
\label{Adprop}
        The number of tilings of an $a$ by $b$ undented Aztec rectangle,
$a<b\leq 2a+1$, $b=2k$, with squares with indices $r_0<r_1<\dots <r_{b-a-1}$
missing from a diagonal of length $a+1$ 
touching the central square is
\[
\frac{2^{k^2-k+a}}{V_{k-1} V_k} 
\biggl(
\prod_{0\leq j<i<2k+1-a} (r_i-r_j ) 
\prod_{\substack{0\leq i<2a-2k+1\\[1pt] 0\leq j<2k-a}} \left| t_i-r_j \right| \biggr)
\determinant{ t_{2i}^j }_0^{a-k} 
\determinant{ t_{2i+1}^j }_0^{a-k-1},
\]
where $t_0<t_1<\dots <t_{2a-2k}$ are the elements of $\{ 0,1,\dotsb a\} -
\{r_0,r_1,\dots ,r_{2k-a-1} \}$.
\end{Prop}

\section{Computing Determinants of Aztec Rectangles: A Special Case}

We can now solve Problem \ref{Propp10} using Proposition \ref{Adprop} 
with $a = 2k-1$, $b=2k$, $r_0=k-1$. The number of tilings is 
$$
\begin{aligned}
\frac{2^{k^2-k+a}}{V_{k-1} V_k}  &
\biggl(\prod_{\substack{0\leq i<2a-2k+1\\0\leq j<2k-a}} \left| t_i-r_j \right| \biggr)
\determinant{ t_{2i}^j }_0^{a-k} 
\determinant{ t_{2i+1}^j }_0^{a-k-1} \\
&= \frac{2^{k^2+k-1}}{V_{k-1} V_k}  
\biggl(\prod_{0\leq i<2k-1} \left| t_i-(k-1) \right| \biggr)
\determinant{ t_{2i}^j }_0^{k-1} 
\determinant{ t_{2i+1}^j }_0^{k-2} \\
&= \frac{2^{k^2+k-1}}{V_{k-1} V_k}  
(k-1)!\, k!\,
\determinant{ t_{2i}^j }_0^{k-1} 
\determinant{ t_{2i+1}^j }_0^{k-2} \\
&= \frac{2^{k^2+k-1}}{V_{k-2} V_{k-1}}  
\prod_{0\leq j_0 < j_1 < k}
(t_{2j_1}-t_{2j_0})
\prod_{0\leq j_0 < j_1 < k-1}
(t_{2j_1+1}-t_{2j_0+1}) \\
&= \frac{2^{k^2+k-1}}{V_{k-2} V_{k-1}}  
\prod_{0\leq j_0 < j_1 < k}
(t_{2j_1}-t_{2j_0})
\prod_{0\leq j_0 < j_1 < k-1}
(t_{2j_1+1}-t_{2j_0+1}) .
\end{aligned}
$$

For $k=2q$, we have
$$
\begin{aligned}
\prod_{0\leq j_0 < j_1 < k}
(t_{2j_1}&-t_{2j_0})\\
&=\prod_{0\leq j_0 < j_1 < q} (2 j_1-2 j_0)  
\prod_{0\leq j_0 < j_1 < q} \bigl((2q+1+2 j_1)-(2q+1+2 j_0)\bigr) \\
&\qquad\qquad\times\prod_{0\leq j_0,j_1 < q} 
((2q+1+2 j_1)-(2 j_0)) \\
&= \bigl(2^{q-1} 4^{q-2} \cdots (2q-2)\bigr)^2 \\
&\qquad\qquad\times3 \cdot 5^2 \cdots (2q-1)^{q-1} (2q+1)^q (2q+3)^{q-1} \cdots
(4q-1)
\end{aligned}
$$
and
\[
\begin{aligned}
\prod_{0\leq j_0 < j_1 < k-1}
(t_{2j_1+1}-t_{2j_0+1})\hskip -1in&\\
&= \prod_{0\leq j_0 < j_1 < q-1} ((2 j_1+1)-(2 j_0+1)) 
\prod_{0\leq j_0 < j_1 < q} ((2q+2 j_1)-(2q+2 j_0)) \\
&\qquad\qquad\times\prod_{\substack{0\leq j_0<q-1\\[1 pt] 0\leq j_1 < q}}
((2q+2 j_1)-(2 j_0+1)) \\
&= (2^{q-2} 4^{q-3} \cdots (2q-4))   (2^{q-1} 4^{q-2} \cdots (2q-2))\\
&\qquad\qquad \times3\cdot 5^2\cdots (2q-1)^{q-1} (2q+1)^{q-1}\cdots (4q-3).
\end{aligned}
\]

For $k=2q+1$, we have
$$
\begin{aligned}
\prod_{0\leq j_0 < j_1 < k}
(t_{2j_1}-&t_{2j_0}) \\
&=\prod_{0\leq j_0 < j_1 < q} (2 j_1-2 j_0)  
\prod_{0\leq j_0 < j_1 < q+1} ((2q+1+2 j_1)-(2q+1+2 j_0)) \\
&\qquad\qquad\times\prod_{\substack{ 0\leq j_0 < q \\[1pt] 0\leq j_1 < q+1}}
((2q+1+2 j_1)-(2 j_0)) \\
&=
(2^{q-1} 4^{q-2} \cdots (2q-2)) 
(2^q 4^{q-1} \cdots (2q))  \\
&\qquad\qquad\times 3 \cdot 5^2 \cdots (2q-1)^{q-1} (2q+1)^q (2q+3)^q \cdots (4q+1)
\end{aligned}
$$
and
$$
\begin{aligned}
\prod_{0\leq j_0 < j_1 < k-1}
(t_{2j_1+1}-t_{2j_0+1})\hskip -1in& \\
&= \prod_{0\leq j_0 < j_1 < q} ((2 j_1+1)-(2 j_0+1)) 
\prod_{0\leq j_0 < j_1 < q} ((2q+2+2 j_1)-(2q+2+2 j_0)) \\
&\qquad\qquad\times\prod_{0\leq j_0,j1<q}
((2q+2+2 j_1)-(2 j_0+1)) \\
&= 
(2^{q-1} 4^{q-2} \cdots (2q-2))^2 \\
&\qquad\qquad\times 3\cdot 5^2\cdots (2q-1)^{q-1} (2q+1)^q (2q+3)^{q-1}\cdots (4q-1)
\end{aligned}
$$

Therefore, for $k=2q$ the number of tilings is
\begin{multline*} 
\frac{2^{(2q)^2+2q-1}}{V_{2q-2} V_{2q-1}} 
2^{4q-5} 4^{4q-9} \cdots (2q-2)^3 \\
\times 3^2 5^4 \cdots (2q-1)^{2q-2} (2q+1)^{2q-1} (2q+3)^{2q-3} \cdots 
(4q-3)^3 (4q-1),
\end{multline*}
and for $k=2q+1$ the number of tilings is
\begin{multline*}
\frac{2^{(2q+1)^2+(2q+1)-1}}{V_{2q-1} V_{2q}}
2^{4q-3} 4^{4q-7} \cdots (2q-2)^5 (2q) \\
\times3^2 5^4 \cdots (2q-1)^{2q-2} (2q+1)^{2q} (2q+3)^{2q-1} \cdots (4q-1)^3 
(4q+1).
\end{multline*}

\section{Computing Determinants: Hexagons}

In this section we solve Propp's Problem 1, and more generally, we count tilings of a
$(2m-1,2n,2m-1)$ or $(2m,2n-1,2m)$ hexagon with a vertical lozenge at the
center. (A
$(k,q,k)$ hexagon has a central vertical lozenge if and only if $k+q$ is odd.)

\begin{Lem}\label{Lem:Vert}
The number of tilings of a $(2m-1,2n,2m-1)$ hexagon with a vertical
lozenge in the center is
\[\frac{1}{V_{2m-2}^2}\determinant{(1+(-1)^{i+j})S_{m+n-1}^{i+j}}_1^{2m-2}.\]

The number of tilings of a $(2m,2n-1,2m)$ hexagon with a vertical lozenge
in the center is 
\[\frac{1}{V_{2m-1}^2}\determinant{(1+(-1)^{i+j})S_{m+n-1}^{i+j}}_1^{2m-1}.\]
\end{Lem}

\begin{proof}
By Proposition \ref{Prop:L}, the number of tilings of a $(2m-1,2n,2m-1)$ hexagon is 
\[\frac{1}{V_{2m-2}^2}  \biggdet{ 
  \sum_{l=0}^{2m+2n-2} l^{i+j}}_0^{2m-2}.\]
By Lemma \ref{Lem:Shift}, this determinant is equal to
\begin{align*}\biggdet{ \sum_{l=-m-n+1}^{m+n-1} l^{i+j}}_0^{2m-2}
&=\determinant{(1+(-1)^{i+j})S_{m+n-1}^{i+j}+\delta_{i+j}}_0^{2m-2}\\
&=\determinant{(1+(-1)^{i+j})S_{m+n-1}^{i+j}}_0^{2m-2}
+\determinant{(1+(-1)^{i+j})S_{m+n-1}^{i+j}}_1^{2m-2},
\end{align*}
where $\delta_k$ is 1 if $k=0$ and is 0 otherwise.

It also follows from Proposition \ref{Prop:L} that the number of tilings of a
$(2m-1,2n,2m-1)$ hexagon that do not have a vertical lozenge in the center is 
\[\frac{1}{V_{2m-2}^2}
\biggdet{ \sum_{\substack{0\le l \le 2m+2n-2\\l\ne m+n-1}} l^{i+j}}_0^{2m-2}.
\]
By Lemma \ref{Lem:Shift} this determinant is equal to 
\[\biggdet{ \sum_{\substack{-m-n+1\le l \le m+n-1\\l\ne 0}} l^{i+j}}_0^{2m-2}
=\determinant{(1+(-1)^{i+j})S_{m+n-1}^{i+j}}_0^{2m-2}.\]
We find the number of tilings that {\it do} have a vertical lozenge in the center by
subtracting from the total number of tilings the number of tilings that do not have a
lozenge in the center.

The formula for $(2m,2n-1,2m)$ hexagons is derived similarly.
\end{proof}

As a first step in evaluating the determinants in Lemma \ref{Lem:Vert}, 
we evaluate the determinant $\determinant{ S_p^{i+j} }_0^{k-1}$. 
It is interesting to note that this determinant was evaluated by
Zavrotsky \cite{Z85} in the course of his research on
minimum square sums, and we follow his proof.

\begin{Lem}\label{Lem:Z}
\[
\begin{aligned}
\determinant{ S_p^{i+j} }_0^{k-1} 
  &=\frac{V_{k-1}^4}{V_{2k-1}}(p-k+1)
  \cdots (p-1)^{k-1}p^k (p+1)^{k-1}\cdots (p+k-1)\\
  &= \frac{V_{p+k-1}V_{p-k-1}V_{k-1}^4}{V_{p-1}^2 V_{2k-1}},
\end{aligned}
\]
 where $S_p^i = 1^i+2^i+\dots+ p^i$ and $V_p=1!\, 2!\cdots p!$.
\end{Lem}
\begin{proof}[Proof {\rm(\/}Zavrotsky \cite{Z85}{\rm)\/}]
If $p$ is a positive integer, we can express the
matrix $(S_p^{i+j} )_0^{k-1}$ as the product of a $k$ by $n$ matrix and  an $p$ by $k$  
matrix,  as in the proof of Proposition \ref{Prop:L}. Since the rank of an
$p$ by $k$ matrix is at most $p$, the rank of the matrix
$(S_p^{i+j} )_0^{k-1}$ is at most $p$. Moreover, this holds also for $p=0$.

Now let $ (a_{i,j} (\lambda ))_0^{k-1}$ be a
matrix whose entries  are polynomials in $\lambda$. 
It is known \cite[p.~17]{FD47} that if,
for some value $\lambda_0$ of $\lambda$, the matrix
$( a_{i,j} (\lambda_0 ))_0^{k-1}$ has rank at most $m$, where $m\le k$,
then $\determinant{a_{i,j} (\lambda )}_0^{k-1}$ is divisible by 
$(\lambda - \lambda_0)^{k-m}$ as a polynomial in $\lambda $.

By Lemma \ref{Lem:0}, there is a polynomial $S_\lambda^i$ in $\lambda$ whose value at
$\lambda=p$ is $S_p^i$. Since the rank of $(S_m^{i+j} )_0^{k-1}$ is at most $m$,
$\determinant{ S_\lambda^{i+j} }_0^{k-1} $ is divisible by $(\lambda -m)^{k-m}$ for 
$0\le m\le k$.

Since $S_{-m}^i = (-1)^{i+1} S_{m-1}^i$ when $i>0$, it follows that the rank of 
$(S_{-m}^{i+j})_0^{k-1}$ is at most one more than the rank of $(S_{m-1}^{i+j})_0^{k-1}$;
i.e., at most $m$. Thus $\determinant{ S_\lambda^{i+j} }_0^{k-1} $ is divisible by 
$(\lambda + m)^{k-m}$ for $1\le m\le k$.

Since $S_\lambda^i$ is a polynomial in $\lambda$ of degree $i+1$, 
$\determinant{ S_\lambda^{i+j} }_0^{k-1}$ is a polynomial in $\lambda$ of degree $k^2$.
Hence 
$\determinant{ S_\lambda^{i+j} }_0^{k-1} $ is equal to a constant times
\[(\lambda-k+1)\cdots 
(\lambda-1)^{k-1} \lambda^k (\lambda+1)^{k-1}
(\lambda+2)^{k-2} \cdots (\lambda+k-1).\]  Since
 $S_\lambda^i$ has
leading coefficient $1/(i+1)$, and, by 
\cite[p.~425]{Ar}, the determinant
 $\determinant{ 1/(i+j+1) }_0^{k-1}$
(a Hilbert determinant) is equal to $V_{k-1}^4/V_{2k-1}$,
we may
compare leading coefficients and  that the constant is
$V_{k-1}^4/V_{2k-1}$.
\end{proof}

\begin{Cor}\label{Cor:Hex}
The number of tilings of a $(k,q,k)$ hexagon is 
\[\frac{V_{2k+q-1}V_{q-1}V_{k-1}^2}{V_{k+q-1}^2V_{2k-1}}.\]
In particular, the number of tilings of a $(2m-1,2n,2m-1)$ hexagon is 
\[\frac{V_{4m+2n-3}V_{2n-1}V_{2m-2}^2}{V_{2m+2n-2}^2 V_{4m-3}}\]
and the number of tilings of a $(2m,2n-1,2m)$ hexagon is 
\[\frac{V_{4m+2n-2}V_{2n-2}V_{2m-1}^2}{V_{2m+2n-2}^2 V_{4m-1}}.\]
\end{Cor}
\begin{proof}
By Proposition \ref{Prop:L} and Lemma \ref{Lem:Shift}, the number of tilings of
a $(k,q,k)$ hexagon is 
\[\frac1{V_{k-1}^2}\determinant{S_{k+q}^{i+j}}_0^{k-1}.\] The result then follows from Lemma
\ref{Lem:Z}.
\end{proof}

It is also possible, as shown in \cite{CLP}, 
to derive the formula for the number of
tilings of an $(a,b,c)$ hexagon directly from Lemma \ref{Lem:1}.

Next we prove a general theorem on Hankel determinants that allows us to evaluate the
determinants in Lemma \ref{Lem:Vert}.

\begin{Prop}\label{Prop:J}
Let $\{ a_i \}_{i=0}^{\infty }$ be a sequence, and let  
\[H_s(k) = \determinant{ a_{(i+j+s)/2}
}_0^{k-1} ,\] 
for $k\ge1$, with $H_s(0)=1$, 
where we take  $a_i$ to be 0 if $i$ is not
an integer. Define $\lambda_k$ inductively by
\[H_0(k+1)=\lambda_0^{k+1}\lambda_1^{k}\cdots \lambda_{k},\]
so that $\lambda_0 = H_0(1)=a_0$  and
\[
\lambda_k = \frac{H_0(k-1) H_0(k+1)}{H_0(k)^2} \] for $k\geq 1$. Then
\begin{equation}
H_2(k)=\lambda_0^{-1}H_0(k+1) \sum_{j=0}^{\lfloor k/2\rfloor} \prod_{i=1}^j
\frac{\lambda_{2i-1} }{\lambda_{2i} } 
\label{Eq:Ho}
\end{equation}
\end{Prop}

\begin{proof}
Define $M_r(k)$ by 
\[M_r(k) = \determinant{ a_{i+j+r} }_0^{k-1}.\] It is easy to see that
\begin{equation}H_{2r}(k)=M_r(\lceil k/2\rceil)M_{r+1}(\lfloor k/2\rfloor).
\label{Eq:M}
\end{equation}
Then
\[\frac{H_2(2m+1)}{H_2(2m)}=\frac{M_1(m+1)}{M_1(m)}
=\frac{H_0(2m+2)}{H_0(2m+1)}.\]
Thus it suffices to prove (\ref{Eq:Ho}) for $k=2m$.

Since (\ref{Eq:Ho}) holds for $k=0$, to prove it for even $k$ 
we need only show that for $m\ge 1$,
\[
\frac{H_2(2m)}{H_0(2m+1)}-\frac{H_2(2m-2)}{H_0(2m-1)}
=\lambda_0^{-1}\prod_{i=1}^m\frac{\lambda_{2i-1}}{\lambda_{2i}}.
\]
Using (\ref{Eq:M}), we may write the identity to be proved as
\begin{equation}
\frac{M_2(m)}{M_0(m+1)}-\frac{M_2(m-1)}{M_0(m)}
=\frac{\lambda_1\lambda_3\cdots\lambda_{2m-1}}{\lambda_0\lambda_2\cdots\lambda_{2m}}.
\label{Eq:M1}
\end{equation}
To prove (\ref{Eq:M1}), we use 
Jacobi's identity \cite[pp.~594--595]{H},
\[
(M_1(m))^2 - M_{0}(m)M_{2}(m) + M_{0}(m+1)M_{2}(m-1) = 0.
\]
Dividing both sides by $M_0(m)M_0(m+1)$, we may rewrite Jacobi's identity as
\begin{equation}
\frac{M_2(m)}{M_0(m+1)}-\frac{M_2(m-1)}{M_0(m)}
=\frac{M_1(m)^2}{M_0(m)M_0(m+1)}.
\label{Eq:J}
\end{equation}
To complete the proof we need to express the right side of (\ref{Eq:J}) in terms of 
the $\lambda_i$.

We have
\begin{align*}
\frac{M_0(m)}{M_0(m-1)}&=\frac{H_0(2m-1)}{H_0(2m-2)}
=\lambda_{2m-2}\frac{H_0(2m-2)}{H_0(2m-3)}\\
&=\lambda_{2m-2}\lambda_{2m-3}\frac{H_0(2m-3)}{H_0(2m-4)}
=\lambda_{2m-2}\lambda_{2m-3}\frac{M_0(m-1)}{M_0(m-2)}.
\end{align*}
Since $M_0(0)=1$ and $M_0(1)=\lambda_0$, this gives
\[\frac{M_0(m)}{M_0(m-1)}=\lambda_0\lambda_1\cdots\lambda_{2m-2},\]
and thus 
\[M_0(m)=\lambda_0^m(\lambda_1 \lambda_2)^{m-1}\cdots (
\lambda_{2m-3} \lambda_{2m-2})
\]

Similarly, we can show that
\[M_1(m) = 
(\lambda_0 \lambda_1)^m\dotsb 
(\lambda_{2m-4} \lambda_{2m-3})^2 (\lambda_{2m-2} \lambda_{2m-1}).
\]

Making these substitutions in the right side of (\ref{Eq:J}) yields (\ref{Eq:M1}),
completing the
proof.
\end{proof}

\begin{note} There is a simple combinatorial proof of Proposition \ref{Prop:J} in which the
determinant is interpreted as counting nonintersecting paths; see Viennot \cite[Chapter~IV]{V}.
\end{note}

We now apply Proposition \ref{Prop:J} to evaluate the determinant
$\determinant{(1+(-1)^{i+j})S_p^{i+j}}_1^k$.

\begin{Prop}
\label{Prop:det}
The determinant $\determinant{(1+(-1)^{i+j})S_p^{i+j}}_1^k$ is equal to 
\[\frac1{2p+1} 
\frac{V_{2p+k+1}V_{2p-k-1}V_{k}^4}{V_{2p}^2V_{2k+1}}
\sum_{j=0}^{\lfloor{k/2}\rfloor}\frac{(\frac12)_j^2\,(\frac
54)_j\,\hskip 5pt(-p)_j\,\hskip3.7pt(p+1)_j}
{(1)_j^2\,(\frac14)_j\,(\frac32+p)_j\,(\frac12-p)_j},
\]
where $(a)_j=a(a+1)\cdots(a+j-1)$.
\end{Prop}

\begin{proof}
Let us set $a_i=2S_p^{2i}+\delta_i$.
Then 
\[a_{i/2}
=\begin{cases} 2S_p^i+\delta_{i/2},& \text{if $i$ is even}\\ 
             0,&\text{if $i$ is odd}\end{cases}
=(1+(-1)^i)S_p^i+\delta_i.
\]

With the notation of Proposition \ref{Prop:J}, the determinant to be evaluated is 
\[\determinant{a_{(i+j)/2}}_1^k=\determinant{a_{(i+j+2)/2}}_0^{k-1}=H_2(k).\] Thus
by Proposition \ref{Prop:J} we can express the value of this determinant in terms
of the values of 
the corresponding determinants $H_0(k)$.

We have
\[(1+(-1)^i)S_p^i+\delta_i
=\sum_{l=-p}^p l^i,\]
so by Lemma \ref{Lem:Shift}, the determinant
$H_0(k)=\determinant{a_{(i+j)/2}}_0^{k-1}$ is equal to
$\determinant{S_{2p+1}^{i+j}}_0^{k-1}$. This determinant may be evaluated by Lemma \ref{Lem:Z}, which 
gives
 \[H_0(k)=\frac{V_{2p+k}V_{2p-k}V_{k-1}^4}{V_{2p}^2V_{2k-1}}.\]
Therefore, with $\lambda_k$ as in Proposition \ref{Prop:J}, we have
$\lambda_0=a_0=2S_p^0=2p+1$, and for $k>0$,
\[\lambda_k=\frac{k^2}4 \frac{(2p+k+1)(2p-k+1)}{(2k-1)(2k+1)}.\]
Thus by Proposition \ref{Prop:J}, we have 
\begin{align*}H_2(k)&=
\lambda_0^{-1}H_0(k+1) \sum_{j=0}^{\lfloor k/2\rfloor} \prod_{i=1}^j
\frac{\lambda_{2i-1} }{\lambda_{2i} } \\
&=\frac1{2p+1} 
\frac{V_{2p+k+1}V_{2p-k-1}V_{k}^4}{V_{2p}^2V_{2k+1}}
\sum_{j=0}^{\lfloor{k/2}\rfloor}\frac{(\frac12)_j^2\,(\frac
54)_j\,\hskip 5pt(-p)_j\,\hskip3.7pt(p+1)_j}
{(1)_j^2\,(\frac14)_j\,(\frac32+p)_j\,(\frac12-p)_j},
\end{align*}
where $(a)_j=a(a+1)\cdots(a+j-1)$.
\end{proof}

We can now combine Lemma \ref{Lem:Vert} with the determinant evaluation
of Proposition \ref{Prop:det} to count tilings of hexagons with a
vertical lozenge in the center:

\begin{Thm} \label{Thm:Hex}
The number of tilings of a $(2m-1,2n,2m-1)$ hexagon with a vertical
lozenge in the center is
\[\frac{V_{4m+2n-3}V_{2n-1}V_{2m-2}^2}
{(2m+2n-1)V_{2m+2n-2}^2 V_{4m-3}}
\sum_{j=0}^{m-1}\frac{(\frac12)_j^2\,(\frac
54)_j\,(1-m-n)_j\,\hskip9pt(m+n)_j\hskip9pt}
{(1)_j^2\,(\frac14)_j\,(\frac12+m+n)_j\,(\frac32-m-n)_j},
\]
and the number of tilings of a $(2m,2n-1,2m)$ hexagon with a vertical
lozenge in the center is  
\[\frac{V_{4m+2n-2}V_{2n-2}V_{2m-1}^2}
{(2m+2n-1)V_{2m+2n-2}^2 V_{4m-1}}
\sum_{j=0}^{m-1}\frac{(\frac12)_j^2\,(\frac
54)_j\,(1-m-n)_j\,\hskip9pt(m+n)_j\hskip9pt}
{(1)_j^2\,(\frac14)_j\,(\frac12+m+n)_j\,(\frac32-m-n)_j}
.\]
\end{Thm}

To finish the solution of  Propp's Problem
1, we need only evaluate the sum in Theorem \ref{Thm:Hex} in the case 
$m=n$.  To do this we use the 
Wilf-Zeilberger (WZ) method \cite{WZ90}.

\begin{Lem} \label{Lem:WZ}
\[\sum_{i=0}^{n-1}\frac{(\frac12)_i^2\,(\frac
54)_i\,(1-2n)_i\,\hskip9.3pt(2n)_i\hskip9.3pt}
{(1)_i^2\,(\frac14)_i\,(\frac12+2n)_i\,(\frac32-2n)_i}
 = \frac{4n-1}{3}.\] 
\end{Lem}
\begin{proof}
Let \[Q(n,i)=\frac 1{4n-1}\frac{(\frac12)_i^2\,(\frac
54)_i\,(1-2n)_i\,\hskip9.3pt(2n)_i\hskip9.3pt}
{(1)_i^2\,(\frac14)_i\,(\frac12+2n)_i\,(\frac32-2n)_i}.\]
We want to prove that 
\[\sum_{i=0}^{n-1}Q(n,i)={\frac{1}{3}} .\]
Since this identity is clearly true for $n=1$, it is sufficient to prove
that
\[\sum_{i=0}^{n}Q(n+1,i) -\sum_{i=0}^{n-1}Q(n,i)=0\] for $n>1$.

To apply the WZ method, we must first find a function $U(n,i)$ such that 
\begin{equation} \label{Eq:QU}
U(n,i+1)-U(n,i)=Q(n+1,i)-Q(n,i).
\end{equation}
With the help of Maple, we find that if we set
\[U(n,i)= {\frac {{i}^{2}\left (2i+1-4n\right
)\left (1+4n\right )\left ( 8{n}^{2}+4n-2{i}^{2}+i+1\right )}{\left
(4i+1\right )\left (2 i+1+4n\right )\left (i-2n\right )\left
(i-1-2n\right )\left (2 n+1\right )n}}Q(n,i)\] 
then (\ref{Eq:QU}) is satisfied.
(Once we have this formula for $U(n,i)$, the verification of (\ref{Eq:QU})
is straightforward.)

Next, we sum identity (\ref{Eq:QU}) on $i$ from 0 to $n-1$ and add $Q(n+1,n)$ to
both sides. The left side telescopes, and we get
\begin{equation}\label{Eq:2}
Q(n+1,n)+U(n,n)-U(n,0)=\sum_{i=0}^{n}Q(n+1,i)
-\sum_{i=0}^{n-1}Q(n,i).
\end{equation}
But $U(n,0)=0$ and we can easily check that $Q(n+1,n)+U(n,n)=0$. Thus the left
side of $(\ref{Eq:2})$ is 0, hence so is the right side.
\end{proof}

\begin{note} The sum in Lemma \ref{Lem:WZ} is a partial sum of a special case of Dougall's
very-well-poised $_5F_4(1)$ sum
\cite[p.~25, eq.~(3)]{B}: if the
upper limit of summation were $2n -1$ instead of $n-1$, we would have a
special case of Dougall's theorem. It is interesting to note that in Ciucu and
Krattenthaler's solution of Propp's Problem 1, they used an
analogous evaluation of a partial sum of the Pfaff-Saalsch\"utz theorem
\cite[eq. (2.1)]{CK}.
\end{note}

We can now finish our solution to Propp's Problem 1:

\begin{Thm}
In a $(2n-1,2n,2n-1)$ or $(2n,2n-1,2n)$ hexagon, the two central
triangles are covered  by a lozenge in exactly one-third of the tilings.
\end{Thm}
\begin{proof}
We compare the result of setting $m=n$ in Corollary \ref{Cor:Hex} with the result of
setting $m=n$ in Theorem \ref{Thm:Hex} and evaluating the sum by Lemma \ref{Lem:WZ}.
\end{proof}

\section{Computing More Determinants}

By  Proposition \ref{Prop:L2}, 
evaluating the Hankel determinant
$\determinant{(1+(-1)^i)S_n^{i+j+1}}_0^{k-1}$ will solve Propp's Problem 2. To
do this, we use the close connection between Hankel determinants and
continued fractions that was implicit in our proof of  Proposition
\ref{Prop:J}. The following lemma is equivalent to \cite[Thm. 7.2]{JT80}.
\begin{Lem}\label{Lem:cf1}
Let $\{ a_i \}_{i=0}^{\infty }$ be a sequence, and 
suppose that the generating function for the $a_i$ has the continued fraction
\[\sum_{i=0}^\infty a_i x^i=\cfrac{\lambda_0}{1-\cfrac{
\lambda_1x}{1-\cfrac{\lambda_2x}{1-\cfrac{\lambda_3x}{ 1-\cdots}}}}
\]
Then
\[
\determinant{ a_{(i+j)/2}}_0^{k-1}=\lambda_0^k\lambda_1^{k-1}\cdots \lambda_{k-1},
\]
where we take  $a_r$ to be 0 if $r$ is not
an integer. \qed
\end{Lem}

By Lemma \ref{Lem:cf1}, if we can find the continued fraction for
\[\sum_{j=0}^\infty (1+(-1)^j)S^{j+1}_n x^{j/2}=\sum_{i=0}^\infty
2S_n^{2i+1}x^{i}, 
\]
then we can evaluate the corresponding Hankel determinant.

The continued fraction in question is given by the following formula, which
we prove in the next section.

\begin{Prop}\label{Prop:cf2}
\[\sum_{i=0}^\infty 2S_n^{2i+1}x^{i}
=\cfrac{\mu_0 }{1-\cfrac{\mu_1x}{1-\cfrac{\mu_2x}{1-\cdots}}},
\]
where
\begin{align*}
\mu_0&=n(n+1)\\
\mu_{2i}&=\frac{i}{4i+2}(n+i+1)(n-i),\quad i\ge 1\\
\mu_{2i+1}&=\frac{i+1}{4i+2}(n+i+1)(n-i).
\end{align*}
\end{Prop}

Now from Proposition \ref{Prop:L2}, Lemma \ref{Lem:cf1}, and Proposition \ref{Prop:cf2}, we
obtain the solution to Propp's Problem 2: 

\begin{Thm}\label{Thm:notri}
The number of tilings of a $(k,2n+1-k,k,k+1,2n-k,k+1)$ hexagon without
the central triangle is
\begin{gather*}\frac{(n-q) (n-q+1)^5\cdots n^{4q+1} (n+1)^{4q+1}\cdots (n+q+1) }
{2^{(2q)(2q+1)} 3^{8(q-1)+2} 5^{8(q-2)+2} \cdots (2q+1)^2}\quad\text{for $k=2q+1$,}\\
\frac{(n-q+1)^3\cdots n^{4q-1} (n+1)^{4q-1}\cdots (n+q)^3}
{2^{(2q-1)(2q)} 3^{8(q-1)-2} 5^{8(q-2)-2} \cdots (2q-1)^6}\quad\text{for $k=2q$.}
\end{gather*}
\end{Thm}

\section{Proof of the continued fraction}

In this section we prove the continued fraction of Proposition \ref{Prop:cf2}.

The exponential generating function for $(1+(-1)^{j-1})S_n^j$ is, by Lemma
\ref{Lem:0},
\begin{align*}
\sum_{j=0}^\infty (1+(-1)^{j-1}) S_n^j \frac{x^j}{j!}
  &=\frac{e^x(e^{nx}-1)}{e^x-1} -\frac{e^{-x}(e^{-nx}-1)}{e^{-x}-1}\\
  &=\frac{(e^{nx}-1)(e^x-e^{-nx})}{e^x-1}\\
  &=2\frac{\sinh\frac n2x \sinh\frac{n+1}2x}{\sinh\frac x2}.
\end{align*}
Now let $L$ be the linear operator on formal power series defined by
\[L\biggl(\sum_{i=0}^\infty u_i\frac{x^i}{i!}\biggr) =\sum_{i=0}^\infty u_ix^i.\]
We note that $L(f(x))$ has the ``formal" integral representation 
\[L(f(x))=\frac 1x \int_0^\infty f(t) e^{-t/x}\,dt,\]
obtained by performing the integration term by term. If $F(x)=L(f(x))$, then this
formula may be written as a Laplace transform
\[F(1/z) = z\int_0^\infty f(t) e^{-tz}\,dt,\]
and this is the form in which it is most often seen in the literature on continued
fractions.

The continued fraction we need is given by the following formula, in which $n$ need not be an
integer. The case in which $n$ is a nonnegative integer is clearly equivalent to Lemma
\ref{Prop:cf2}.

\begin{Lem}\label{Lem:cf3}
\[L\left(\frac{\sinh\frac n2x \sinh\frac{n+1}2x}{\sinh\frac x2}\right)
 = \cfrac{\binom{n+1}{2}x}{1-\cfrac{\mu_1x^2}{1-\cfrac{\mu_2x^2}{1-\cdots}}},
\]
where
\begin{align*}
\mu_{2i}&=\frac{i}{4i+2}(n+i+1)(n-i)\\
\mu_{2i+1}&=\frac{i+1}{4i+2}(n+i+1)(n-i).
\end{align*}
\end{Lem}

Lemma \ref{Lem:cf3} is one of several continued fractions for this function given by Lange
\cite[pp.~259--260]{La}. (A closely related continued fraction for the same function was given by
Stieltjes \cite{St}.) For completeness, we give here a self-contained proof:

\let\a=n
\def\floor#1{\left\lfloor#1\right\rfloor}
\def\F#1#2#3#4{{}_2F_1\left(\left.\genfrac{}{}{0pt}{}{#1,\,#2}{#3}\,\right|
#4\right)}

\begin{Lem}\label{cf}
Let $f_0, f_1, f_2,\dots$ be formal power series in $x$  with nonzero constant terms, and let
$c_1,c_2,\dots$ be constants such that for each
$k\ge1$,
\begin{equation}\label{recur}
f_k-f_{k-1}=c_k x^2  f_{k+1}.
\end{equation}
Then for each $m\ge1$,
\[\frac {f_m}{f_{m-1}}=
  \cfrac1{1-\cfrac{c_m x^2}{1-\cfrac{c_{m+1}x^2}{1-\cdots}}}
\]
\end{Lem}
\begin{proof}
Equation \eqref{recur} is equivalent to 
\[\frac{f_k}{f_{k-1}}=\cfrac{1}{1-c_k x^2\cfrac{f_{k+1}}{f_k}}.\]
Iterating this formula gives
\[\frac{f_m}{f_{m-1}}=\cfrac{1}{1-
\cfrac {c_m x^2}{1-
\cfrac{c_{m+1} x^2}
{\cfrac{\ddots}{1-c_{m+n} x^2 \cfrac{f_{m+n+1}}{f_{m+n}}}}}}
\]
Taking the limit as $n\to\infty$ yields the lemma. 
\end{proof}

\begin{proof}[Proof of Lemma \ref{Lem:cf3}]

Let $E=L^{-1}$, so that 
\[E\biggl(\sum_{j=0}^\infty a_j x^j\biggr)=\sum_{j=0}^\infty a_j \frac{x^j}{j!},\]
and 
suppose that with $f_k$ as in Lemma \ref{cf}, $g_k=E(x^kf_k)$.
Multiplying the recurrence \eqref{recur} by $x^{k-2}$, and using the fact that if
$u(x)$ is divisible by $x$ then 
\[E\left(\frac{u(x)}x\right)=\frac{d\ }{dx}E\bigl(u(x)\bigr),\]
we find that (\ref{recur}) is equivalent to
\begin{equation}\label{diffeq}
\frac{dg_{k}}{dx}-g_{k-1}=c_{k}g_{k+1}.
\end{equation}

We now consider the case of Lemma \ref{cf} in which 
\begin{align*}
c_{2i}&=\frac i {4i+2}(\a+i+1)(\a-i)\\
c_{2i+1}&=\frac{i+1}{4i+2}(\a+i+1)(\a-i)\\
\end{align*}

We shall express a solution of  recurrence \eqref{diffeq} in terms 
of the hypergeometric series, defined by
\[ \F abcz = \sum_{n=0}^\infty \frac{(a)_n(b)_n}{n!\,(c)_n}z^n,\]
where $(u)_n=u(u+1)\cdots(u+n-1)$.

We claim that a solution of  recurrence \eqref{diffeq} is
\begin{multline}\label{sol} \qquad g_k=\frac{(e^x-1)^k}{k!}e^{-\a x}
\F{\floor{\frac{k+1}2}}
{\floor{\frac{k+1}2} -\a}
{k+1}{1-e^x}\\
\times\F{\floor{\frac{k}2}+1}
{\floor{\frac{k}2} -\a}
{k+1}{1-e^x}.\qquad
\end{multline}
The verification that $g_i$ defined by (\ref{sol}) really does satisfy 
\eqref{diffeq} is a straightforward, but tedious, computation using the formula
\[ \frac{d\ }{dz}\F abcz=\frac{ab}c\F{a+1}{b+1}{c+1}z,\]
together with the contiguous relations for the hypergeometric series
\cite[p.~558]{AS}. (This computation was done with
the help of Maple.) 

It is not hard to show that $g_0=1$. We now evaluate 
\[g_1=(e^x-1) e^{-\a x}\F1{1-\a}2{1-e^x}
\F1{-\a}2{1-e^x}
\]
Using the easily verified fact that 
\begin{equation}\label{Eq:2F1}
\F1\beta2z=\frac1{z(\beta-1)}\left(\frac1{(1-z)^{\beta -1}}-1\right),
\end{equation}
we find that
\[g_1=\frac{e^{-\a x}(e^{\a x}-1)(e^{(\a+1)x}-1)}{\a(\a+1)(e^x-1)}
  =\frac{2}{\a(\a+1)}\frac{\sinh{\frac \a2 x}\,\sinh\frac{\a+1}2x}{\sinh\frac x2}.
\]
Thus $f_0=1$, and 
\[f_1=\frac 1xL(g_1)=\frac1{\binom{\a+1}2x}L\left(\frac{\sinh{\frac \a2
x}\,\sinh\frac{\a+1}2x}{\sinh\frac x2}\right).\]
Substituting these values of $f_0$ and $f_1$ into the case $m=1$ of Lemma \ref{cf}, and
multiplying both sides by $\binom{\a+1}2x$, 
completes the proof of Lemma \ref{Lem:cf3}.
\end{proof}

It is clear from the recurrence \eqref{diffeq} and the value of $g_1$ that $g_k$ is 
a rational function of $e^x$ and $e^{\a x}$. Although we won't need
it, we
 can give an explicit formula that expresses $g_k$ in
this form by applying to~\eqref{sol} the formula
\begin{multline*}
 \F{m+1}{\beta}{k+1}z = (-1)^k\frac{k!}{m!\,z^k}
\biggl[
   (-1)^m\frac{(1-z)^{k-m-\beta}}{(1-\beta)_{k-m}}\\
\times\F{-m}{1-\beta}{1-\beta-m+k}{1-z}
- \frac{(1-k)_m}{(1-\beta)_k}\sum_{i=0}^{k-m-1}
\frac{(\beta-k)_i(1-k+m)_i}{i!\,(1-k)_i}z^i\biggr],
\end{multline*}
for $k\ge m$, 
which can be proved by equating coefficients of powers of $z$ on both sides.
(Note that (\ref{Eq:2F1}) is the case $m=0, k=1$.)

\section{Acknowledgments}

We would like to thank Christian Krattenthaler for telling us about the
status of his work on these problems, and James Propp, for his
assistance at various times during the writing of this paper.


\begin{thebibliography}{99}

\bibitem{AS} 
  M. Abramowitz and I. A. Stegun, {\em Handbook of Mathematical Functions,}
Dover, New York, 1972..

        \bibitem{A}
        G. E. Andrews, {\em The Theory of Partitions,}
Encyclopedia of Mathematics and Its
Applications, Vol.~1, Addison-Wesley, Reading, Massachusetts, 1976.

\bibitem{Ar}
J. W. Archbold, {\em Algebra,} Pitman Paperbacks, Bath,1970. 


\bibitem{B}
W. N. Bailey, {\em Generalized Hypergeometric Series,}
Cambridge University Press, London,1935.


       \bibitem{CK}
        M. Ciucu and C. Krattenthaler, ``The number of
rhombus tilings of a symmetric hexagon which contain the central rhombus,''
{\em J. Combin. Theory Ser. A}, to be published.
   
     \bibitem{CLP}
        H. Cohn, M. Larsen, and J. Propp, ``The shape
of a typical boxed plane partition,'' submitted to the {\em New York
Journal of Mathematics}.

        \bibitem{EKLP}
        N. Elkies, G. Kuperberg, M. Larsen, and J. Propp, ``Alternating-Sign
matrices and domino tilings'',  {\em J. Algebraic Combinatorics},
1 (1992),  111--132 and 219--234.

       \bibitem{FK}
        M. Fulmek and C. Krattenthaler, ``The number of rhombus tilings 
of a symmetric hexagon which contain a fixed rhombus on the symmetry axis, I,"
{\em Ann. Combin.} 2 (1998), 19--40.

        \bibitem{FD47}
        R. A. Frazer, W. J. Duncan, and A. R. Collar, {\em Elementary
Matrices and Some Applications to Dynamics and Differential Equations,} Cambridge
University Press, 1947.

\bibitem{G}
F. R. Gantmacher, {\em The Theory of Matrices,} Volume 1, Chelsea, New York,
1959.

        \bibitem{H}
        P. Henrici, {\em Applied and Computational Complex Analysis,} Volume 1,
Wiley, New York, 1974. 

        \bibitem{JT80}
        W. B. Jones and W. J. Thron, {\em Continued Fractions:
Analytic Theory and Applications,} Encyclopedia of Mathematics and Its
Applications, Vol. 11, Addison-Wesley, Reading, Massachusetts, 1980.
        \bibitem{JP}
        W. Jockusch and J. Propp, ``Antisymmetric monotone triangles and
domino tilings of quartered Aztec diamonds," {\em J.
Algebraic Combinatorics,} to be published.

        \bibitem{K}
        P. W. Kasteleyn, ``The statistics of dimers on a lattice,
I. The number of dimer arrangements on a quadratic lattice,''
{\em Physica} 27 (1961), 1209--1225. 

\bibitem{Kn}
        D. E. Knuth, {\em The Art of Computer Programming, Vol. 3: Sorting and
Searching,} Addison-Wesley, Reading, Massachusetts, 1973.

\bibitem{Kr}
C. Krattenthaler, ``Schur function identities 
and the number of perfect matchings of 
holey Aztec rectangles," preprint.
        
        
\bibitem{La}
        L. J. Lange, ``Continued fraction representations for functions
related to the gamma function,'' in {\em Continued Fractions and Orthogonal
Functions: Theory and Applications}, 
ed. S. Clement Cooper and W. J. Thron, M. Dekker, 1994, pp. 233--279.

\bibitem{MRR}
W. H. Mills, D. P. Robbins, and H. Rumsey, Jr., ``Alternating sign matrices and descending
plane partitions," {\em J. Combin. Theory Ser. A} 34 (1983), 340--359.

\bibitem{KO}
        S. Okada and C. Krattenthaler,
``The number of rhombus tilings of a
`punctured' hexagon and the minor summation formula,'' {\em Adv. Appl. Math.}, to be
published

\bibitem{PS}
J. Propp and R. Stanley, ``Domino tilings with  barriers," preprint.

        \bibitem{jP1}
        J. Propp, ``Twenty open problems in enumeration of
matchings,'' online document, 1996:
http://www-math.mit.edu/$\sim$propp/open.ps.gz

        \bibitem{cR79}
        C. Radoux, ``Calcul effectif de certains d\'eterminants
de Hankel,'' {\em Bull. Soc. Math.
Belg. S\'er. B} 31 (1979),  49--55.

\bibitem{St}
T. J. Stieltjes, ``Sur quelques int\'egrales d\'efinies et leur d\'eveloppement
en fractions continues,"{\em Ouvres Compl\`etes,} Vol. 2, P. Noordhoff, Groningen, 1918,
pp. 378--391. Originally published in {\em Quarterly Journal of Math.} 24 (1890), 370--382.

\bibitem{Sy}
J. J. Sylvester, ``On a certain fundamental theorem of determinants," Philosophical Magazine (4) 2 (1851),
142--145.

\bibitem{V}
X. G. Viennot, {\em Une th\'eorie combinatoire des polyn\^omes orthogonaux g\'en\'eraux,}
Conference Notes, Universit\'e du Qu\'ebec \`a Montr\'eal, 1983.

\bibitem{WZ90}
H. S. Wilf and D. Zeilberger,
``Rational functions certify combinatorial identities,''{\em J. Amer.
Math. Soc.} 3 (1990),  147--158.  

\bibitem{Z85}
A. Zavrotsky, ``El Gesseliano'' (Spanish), 
{\em Notas de Matematicas,} no. 73, Universidad de los Andes, Facultad
de Ciencias, Departamento de Matematica, Merida, Venezuela, 1985.

\end{thebibliography}
\end{document}